\documentclass[12pt]{amsart}
\usepackage{amsthm}
\usepackage{epsf}
\usepackage{amssymb}
\usepackage{latexsym}
\usepackage{amsmath}
\newcommand{\moins}{\backslash}
\newcommand{\C}{{\mathbb  C}} 
 
\newcommand{\Z}{{\mathbb  Z}} 
\newcommand{\R}{{\mathbb R}} 
\newcommand{\B}{{\mathbb B}} 
\newcommand{\ca}{{\mathcal A}}
\newcommand{\cb}{{\mathcal B}}

\newtheorem{theorem}{Theorem} %the resolution could also be [subsection]
\newtheorem{definition}[theorem]{Definition}

\newtheorem{corollary}[theorem]{Corollary}
\newtheorem{lemma}[theorem]{Lemma}

\newtheorem{proposition}[theorem]{Proposition}
\newtheorem{remark}[theorem]{Remark}
\def \1{^{-1}}
\def \2{^{-2}}
\def \3{^{-3}}
\def \A{{\mathcal A}}
\def \B{{\mathcal B}}
\def \Cc{{\mathcal C}}
\newcommand\proj{{\mathbb P}^2}
\newcommand\ptwo{{\mathbb P}^2}
\newcommand\comp[1]{\proj\backslash #1}
\newcommand\fg[1]{\pi_1(\comp{#1})}
\newcommand\pline{{\mathbb P}^1}
\newcommand\ok{\to}

\def\C{{\mathbb C}}
\def\P{{\mathbb P}}
\def\R{{\mathbb R}}
\def\Z{{\mathbb Z}}

\long\def\forget#1\forgotten{}
\def \prodl{\prod\limits}

\begin{document}

\renewcommand{\subjclassname}{%
       \textup{2000} Mathematics Subject Classification}
\date{\today}

\address{Meirav Amram, Einstein Mathematics Institute, the Hebrew University, 
Jerusalem, Israel/ Department of Sciences, Holon Academic Institute of Technology, Israel}
\email{ameirav@math.huji.ac.il}

\address{Mina Teicher, Department of Mathematics, Bar-Ilan University, Ramat-Gan, Israel}
\email{teicher@macs.biu.ac.il}

\title {Fundamental groups of some special quadric arrangements}

\author[M.~Amram, M.~Teicher]
{Meirav Amram$^1$ \and Mina Teicher}
\stepcounter{footnote} 
\footnotetext{Partially supported by EU-network HPRN-CT-2009-00099 (EAGER), 
the Emmy Noether Research Institute for Mathematics, the Israel Science Foundation 
grant \#8008/02-3 (Excellency Center ``Group Theoretic Methods in 
the Study of Algebraic Varieties"), the Edmund Landau 
Center for Research in Mathematical Analysis and Related Areas, 
sponsored by the Minerva Foundation (Germany).}

\maketitle
%\stepcounter{footnote}
%\footnotetext
\keywords{{\bf{Keywords:}} Fundamental groups, 
complement of curve, quadric arrangement.}

\subjclass{{\bf{MSC Classification:}} 14H20, 14H30, 14Q05.}

\begin{abstract}
In this work we obtain  presentations of fundamental 
groups of the complements of three families of quadric arrangements in $\proj$. 
The first arrangement is a union of $n$ quadrics, 
which are tangent to each other at two common points. The second arrangement is 
composed of $n$ quadrics  which are tangent to each other at one common point. The 
third arrangement is composed of $n$ quadrics,  $n-1$ of them are tangent to 
the $n$'th one and each one of the $n-1$ quadrics is transversal to the other $n-2$ ones. 
\end{abstract}

\maketitle

\section{Introduction}
The aim of the present article is the computation of the fundamental groups to complements of 
some quadric arrangements in $\ptwo$. Recall that given a quadric-line arrangement in $\ptwo$, 
we are interested in computing the fundamental 
group of its complement. 
The present paper is devoted to the computation of the fundamental 
groups related to three infinite families of quadric arrangements. 
The three types of interesting families of quadric curves in this paper are as 
follows. The first arrangement is a union of $n$ quadrics, 
which are tangent to each other at two common points (Figure \ref{An}). 
The second arrangement is composed of $n$ quadrics  which are tangent to each other at one common point 
(Figure \ref{Bn}). The third arrangement is composed again of $n$ quadrics,  $n-1$ of them are tangent to 
the $n$'th one and each one of the $n-1$ quadrics is transversal to the other $n-2$ ones (Figure \ref{Cn}). 

Some work has been done concerning line arrangements (see e.g. \cite{garber}, 
\cite{orliksalomon}, \cite{suciu}), and other progress has been done also concerning 
quadric-line arrangements (see \cite{amram}, \cite{amramteicher}, and \cite{AmGaTe2}). 

Let $C\subset \proj$ be a plane curve and $* \in \comp{C}$ a base point. 
By abuse of notation, we will call the group 
$\pi_1({\proj\backslash C},*)$ the \textit{fundamental group of C}, 
and we shall frequently omit base points and write $\fg{C}$. 
One is interested in the group $\fg{C}$ mainly for two reasons. First, when the curve 
appears to be a branch curve, then $\fg{C}$ is an important invariant, concerning 
either the branch curve or the surface itself. Secondly, $\fg{C}$ contributes to the study of
the Galois coverings $X\to \proj$ branched along $C$.   
Many interesting surfaces have been constructed as branched Galois coverings 
of the plane. One example is  the arrangement $\A_3$ (shown in Figure \ref{A3} below), 
which has Galois coverings $X\to\proj$ branched along it;  
$X\simeq {\mathbb P}^1\times {\mathbb P}^1$, 
or $X$ is either an abelian surface, a $K3$ surface, or a quotient of the two-ball 
${\mathbb B}_2$ (see \cite{holzapfel}, \cite{uludag1}, \cite{yoshida}). 
Moreover, some line arrangements defined by unitary reflection groups studied in 
\cite{orlik} are related to $\A_3$ via orbifold coverings. 
For example, if ${\mathcal L}$ is the line arrangement given by the equation 
$$
xyz(x+y+z)(x+y-z)(x-y+z)(x-y-z)=0, 
$$
then the image of ${\mathcal L}$ under the branched covering map 
$[x:y:z]\in \proj\to[x^2:y^2:z^2]\in\proj$ is the arrangement $\A_3$,
see \cite{uludag1} for details.
\begin{figure}[h]
\epsfxsize=4cm %width
\epsfysize=3.5cm %height
%\begin{minipage}{\textwidth}
\begin{minipage}{\textwidth}
\begin{center}
\epsfbox{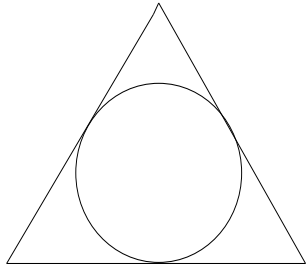}
\end{center}
\end{minipage}
\caption{The arrangement $\A_3$}\label{A3}
\end{figure}

We use the algorithm of Moishezon-Teicher \cite{MoTe} in order to compute the braid monodromy 
of each one of the arrangements. 
Then we use the van Kampen Theorem \cite{vankampen} in order to 
get the presentations of the complements of the arrangements in $\C\P^2$.

This paper is divided into four parts. In Section \ref{notions}, we quote basic definitions and an 
alternative way for computing quadric arrangements, by defining them as birational to line arrangements 
or as their coverings. In Section \ref{results}, we quote the results of this work (Theorems \ref{presAn}, 
\ref{presBn}, \ref{presCn}). In Sections \ref{pres1}, \ref{pres2} and \ref{pres3}, we prove these theorems.

\section{Quadric arrangements related to line arrangements}\label{notions}
\subsection{Meridians}
Let $C\subset X$ be a  curve in a smooth complex surface $X$ and $p\in C$. 
A \textit{meridian} $\mu$ of $C$ at $p$, based on a point $* \in \proj\moins C$
is a loop in $\comp{C}$ obtained by following a path $\omega$ with $\omega(0)=*$ and 
$\omega(1)$ belonging to a small neighborhood of $p$, turning around $C$ in the positive 
sense along the boundary of a small disc $\Delta$, having with $C$ a single intersection at $p$, 
and then turning back to $*$ along $\omega$. 
If $B\subset C$ is an irreducible component, a meridian $\mu_p$ of $C$ at a point  $p\in 
B\moins \mbox{Sing}(C)$ will be called a \textit{meridian} of $B$. 
It is well-known that (homotopy classes of) any two meridians of $B$ are conjugate 
elements in $\pi_1(X\moins C,*)$ (see e.g. \cite[Section 7.5]{lamotke}). 
When $p$ is a singular point of $C$, we have the following result.
\begin{lemma}
Let $p\in C$ be a singular point, $\mu_p$ a meridian of $C$ at $p$, and let 
$\sigma:Y\ok \ptwo$ be the blow-up of $X$ at $p$. Denote by $C$ the proper transform of 
$C$ and by $P$ the exceptional divisor. 
Then $\sigma(\mu_p)$ is a meridian of $P$. In particular, any two meridians of $C$ at $p$ 
are conjugate elements of $\pi_1(X\moins C)\simeq \pi_1(Y\moins(C\cup P))$. 
\end{lemma}

\begin{proof}
The spaces $Y\moins(C\cup P)$ and $X \moins C$ are homeomorphic. 
By definition, $\mu_p=\omega\cdot \partial\Delta\cdot \omega\1$, where $\Delta$ is a disc, 
having an intersection with $C$ at $p$, implying  that the disc $\sigma(\Delta)$ intersects $P$ 
transversally and away from $C$. In other words, the loop $\sigma(\mu)$ 
is a meridian of $P$.  
\end{proof}

The group $\pi_1(X\moins C)$ is an invariant of the pair $(\ptwo, C)$. 
Since meridians are well-defined up to a conjugacy class,  they can be considered as 
supplementary invariants of $\pi_1(X\moins C)$. What follows is a description of  how to 
capture the meridians at singular points of $C$, during the computation 
of the group $\pi_1(X\moins C)$ by the van Kampen Theorem.

\begin{lemma}
Let $C\subset \ptwo$ be a curve, $L_0$ a line in general position with respect to $C$ and 
let $*\in L_0\moins C$ be a base point. Let $p$ be a singular point of  $C$. We assume that 
$L_0$ passes through a sufficiently small ball $V$ around $p$, so that the disc $\Delta:=
L_0\cap V$ meets all branches of $C$ meeting at $p$. Take a path  $\omega$  in $L_0$ 
connecting $*$  to a boundary point $q$ of $\Delta$. Then  $\mu_0:=\omega
\cdot \partial\Delta \cdot \omega\1$ is a homotopic to a meridian of $C$ at $p$. 
\end{lemma}

\begin{proof} 
Let $L_1$ be the line through $*$ and $p$. Consider the projection $\phi:\ptwo\ok\pline$ 
from the point $*$ and  with $\phi(L_0)=a$,  $\phi(L_1)=b$ and take a path $\gamma\subset 
\phi(V)$  with $\gamma(0)=a$ and $\gamma(1)=b$. Put $L_t$ for the fiber above $\gamma(t)$. 
Let $\sigma\subset \partial V$ be a lift of the path $\gamma$ with $\sigma(0)=q$. Assume 
that  $\sigma_\tau$ is the loop $\sigma$ until $\tau$, i.e. with $\sigma_\tau(t):=\sigma
(\tau t)$ ($t\in[0,1]$). Put $\omega_\tau:=\sigma_\tau \cdot \omega$ and define $\mu_\tau:
=\omega_\tau\cdot \partial\Delta_\tau\cdot \omega_\tau\1$, where $\Delta_\tau:=L_\tau\cap V$. 
Then $M(t,\tau):=\mu_\tau(t)$ gives a homotopy between $\mu_0$ and $\mu_1$, and this latter 
loop is obviously a meridian of $C$ at $p$. 
\end{proof}

\subsection{Quadric arrangements birational to line arrangements} 
Assume that $\ca$ is a line arrangement, and let $\psi$ be the involution 
$\psi:[x:y:z] \in \ptwo\ok [1/x:1/y:1/z]\in \ptwo$. Suppose that the lines $X$, $Y$, 
$Z$ are respectively given by the equations $x=0$, $y=0$ and $z=0$. 
If $\ca$ is in general position with respect to $X\cup Y\cup Z$, then  $\psi(\ca)$ is 
an arrangement of smooth quadrics. In addition to those of $\ca$, this arrangement has 
three more singular points where all the irreducible components of $\psi(\ca)$ meet 
transversally. 

In this case, the group $\fg{\psi(\ca)}$ can easily be found in terms of $\fg{\ca}$ as 
follows: Assuming $\ca=\cup_{i=1}^n L_i$, let 
\begin{equation}\label{pres:1}
\fg{\ca}\simeq \langle\mu_1,\dots,\mu_n\,|\, w_1=\cdots=w_m=\mu_1\dots\mu_n=e \rangle 
\end{equation}
be a presentation obtained by an application of van Kampen, where $\mu_i$ is a 
meridian of $L_i$. Put $\ca^\prime :=\ca\cup  X\cup Y\cup Z$. Since $\ca$ is in general 
position with respect to $X\cup Y\cup Z$, one has by \cite{garber}
\begin{equation}\label{pres:2}
\fg{\ca^\prime}\simeq \Biggl\langle 
\begin{array}{l} 
\mu_1,\dots,\mu_n,\\ 
\sigma_1,\sigma_2,\sigma_3 
\end{array} 
\left | 
\begin{array}{l}
[\mu_i,\sigma_j]=[\sigma_j,\sigma_k]=e \;(i\in[1,n],\; j,k\in [1,3])  \\
 w_1=\cdots=w_m=\mu_1\dots\mu_n\sigma_1\sigma_2\sigma_3=e            
\end{array} \right . 
\Biggr\rangle
\end{equation}
where $\sigma_1$, $\sigma_2$, $\sigma_3$ are, respectively, meridians of $X$, $Y$ and $Z$. 
Let $p:=X\cap Y$, $q=Y\cap Z$ and $r:=Z\cap X$.
Then $\sigma_1\sigma_2$ (respectively, $\sigma_2\sigma_3$, $\sigma_3\sigma_1$) is a meridian 
of $\ca^\prime $ at $p$ (respectively, $q$, $r$). 
Hence, the group $ \fg{\psi(\ca)}$ can be obtained by setting $\sigma_1\sigma_2=\sigma_2
\sigma_3=\sigma_3\sigma_1=e$ in the presentation of $\fg{\ca}$. 
But these relations imply $\sigma:=\sigma_1=\sigma_2=\sigma_3$ and $\sigma^2=e$ and 
since by the projective relation one has $\mu_1\dots\mu_n\sigma_1\sigma_2\sigma_3=e$, 
it suffices to  replace this latter relation by $(\mu_1\dots\mu_n)^2=e$. Hence
\begin{equation*}
\fg{\psi(\ca)}\simeq \Biggl\langle 
\begin{array}{l} 
\mu_1,\dots,\mu_n\\ 
\end{array} 
\left | 
\begin{array}{l}
[\mu_i,\mu_1\dots\mu_n]=e \;(i\in[1,n])  \\
 w_1=\cdots=w_m=(\mu_1\dots\mu_n)^2=e            
\end{array} \right . 
\Biggr\rangle.
\end{equation*}
Since $\sigma$ is a central element of this group, this proves the following result.
\begin{theorem} 
For any arrangement of $n$ lines $\ca$, there is an arrangement of $n$ smooth quadrics 
$\cb$ with a central extension 
$$
0\ok \Z/(2) \ok \fg{\cb} \ok \fg{\ca}\ok 0.
$$
\end{theorem}

\subsection{Quadric arrangements as coverings of line arrangements}
Assume that $\ca$ is a line arrangement, and let $\phi$ be the branched covering 
$\phi:[x:y:z] \in \ptwo\ok [x^2:y^2:z^2]\in \ptwo$. Suppose that the lines $X$, $Y$, $Z$ 
are respectively given by the equations $x=0$, $y=0$ and $z=0$. 
If $\ca$ is in general position to $X\cup Y\cup Z$, then  $\phi\1(\ca)$ is an arrangement 
of smooth quadrics. Above any singular point of $\ca$ lie four singular points  of 
$\phi\1(\ca)$ of the same type. 
In this case, the group $\fg{\phi\1(\ca)}$ can easily be found in terms of $\fg{\ca}$ as 
follows: Assuming $\ca=\cup_{i=1}^n L_i$, one has a presentation (\ref{pres:1}).
For the arrangement $\ca^\prime :=\ca\cup  X\cup Y\cup Z$, the presentation (\ref{pres:2}) 
is valid. There is an exact sequence
$$
0\ok \fg{\phi\1(\ca^\prime)}\ok \fg{\ca^\prime}\ok \Z/(2)\oplus\Z/(2)\ok 0. 
$$

The group $\fg{\phi\1(\ca)}$ is the quotient  $\fg{\phi\1(\ca^\prime)}$ by the subgroup 
generated by the meridians of $\phi\1(X)$, $\phi\1(Y)$ and $\phi\1(Z)$.

\section{Statements of results}\label{results}
In this section we give in Theorems \ref{presAn}, \ref{presBn} and \ref{presCn} 
the presentations of the fundamental groups of the three quadric arrangements in $\ptwo$. 
We prove them in the forthcoming sections. For our computations, we need the following definition.
\begin{definition}
A group $G$ is said to be \textit{big} if it contains a non-abelian free subgroup. 
\end{definition}

\bigskip

\subsection{The quadric arrangement $\A_n$}
\begin{figure}[h]
\begin{minipage}{\textwidth}
\begin{center}
\epsfbox{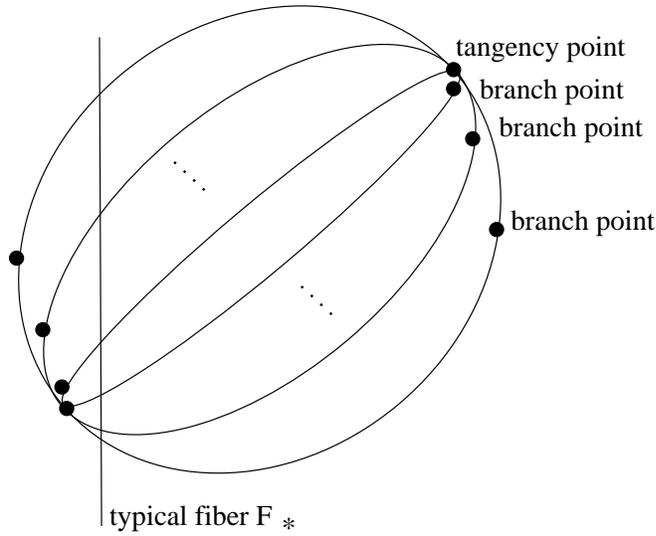}
\end{center}
\end{minipage}
\caption{The arrangement $\A_n$}\label{An}
\end{figure}

Let $\A_n:=Q_1 \cup \cdots \cup Q_n$ be a quadric arrangement, which is a union of $n$ 
quadrics tangent to each other at two common points, see Figure \ref{An}.

\begin{theorem}\label{presAn} 
The fundamental group $\fg{\A_n}$ of $\A_n$ in $\ptwo$ admits the presentation
\begin{equation}\label{equAn}
\fg{\A_n} \simeq \Biggr\langle a_1, a_2, \dots, a_n \;  \left | 
\begin{array}{ll}
(a_1 a_2 \cdots a_n)^2 = e 
\end{array} 
\right. \Biggr\rangle, 
\end{equation}
where $a_1, \dots, a_n$ are meridians of $Q_1, \dots, Q_n$, respectively. 
\end{theorem}

\begin{corollary}\label{bigAn}
The group $\fg{\A_n}$ is abelian for $n=1$ and big for $n \geq 2$. 
\end{corollary}

\bigskip

\subsection{The quadric arrangement $\B_n$}
\begin{figure}[h]
%\epsfxsize=9cm %width
%\epsfysize=5cm %height
%\begin{minipage}{\textwidth}
\begin{minipage}{\textwidth}
\begin{center}
\epsfbox{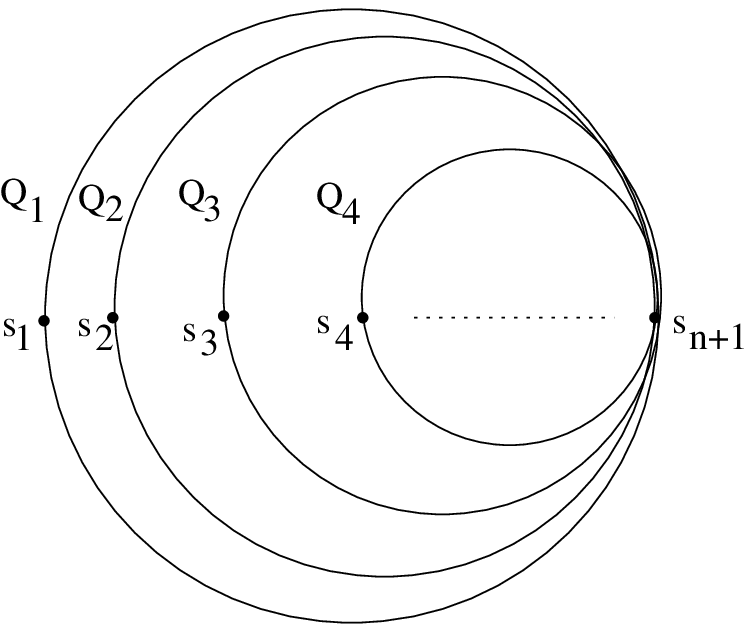}
\end{center}
\end{minipage}
\caption{The arrangement $\B_n$}\label{Bn}
\end{figure}

Let $\B_n:=Q_1 \cup \cdots \cup Q_n$ be a quadric arrangement, composed of $n$ 
quadrics tangent to each other at one common point, see Figure \ref{Bn}. 

\begin{theorem}\label{presBn}
The group $\fg{\B_n}$ admits the presentation
\begin{equation}\label{equBn}
\fg{\B_n} \simeq 
\Biggl\langle 
\begin{array}{l}a_1, a_2, \dots, a_n 
\end{array} 
\left | 
\begin{array}{ll}
 \ (a_1 a_2 \cdots a_n)^2 = e 
\end{array} \right. \Biggr\rangle,
\end{equation}         
where $a_1, \dots, a_n$ are meridians of $Q_1, \dots, Q_n$, respectively.
\end{theorem}

\begin{corollary}\label{bigBn}
The group $\fg{\B_n}$ is abelian for $n=1$ and big for $n \geq 2$. 
\end{corollary}

\bigskip

\subsection{The quadric arrangement $\Cc_n$}
\begin{figure}[h]
\epsfxsize=12cm %width
\epsfysize=8cm %height
\begin{minipage}{\textwidth}
\begin{center}
\epsfbox{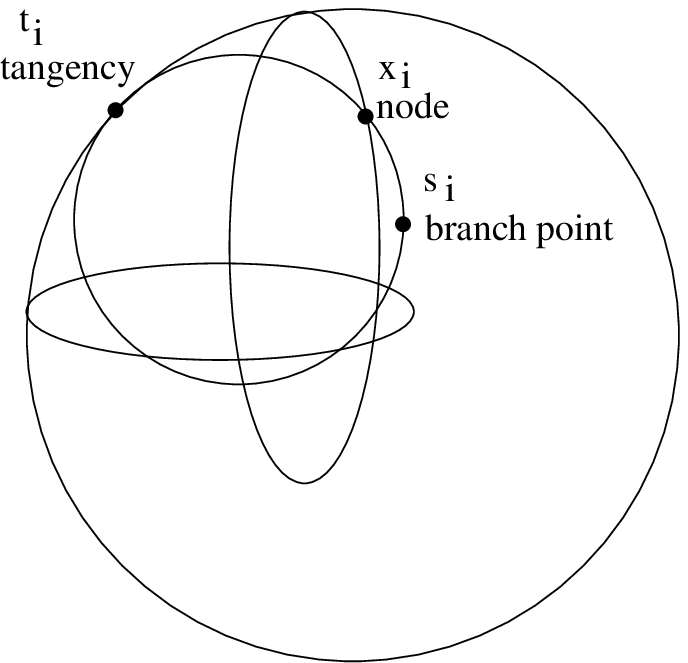}
\end{center}
\end{minipage}
\caption{The arrangement $\Cc_n$}\label{Cn}
\end{figure}

Let $\Cc_n:=Q_1\cup\cdots\cup Q_n$ be a quadric arrangement, which is a union of $n$ 
quadrics,  $n-1$ of them are tangent to the $n$'th one and each one of these $n-1$ 
quadrics is transversal to the other $n-2$, see Figure \ref{Cn}.

\begin{theorem}\label{presCn}
The group $\fg{\Cc_n}$ admits the presentation 
\begin{equation}\label{equCn}
\fg{\Cc_n} \simeq \Biggr\langle a_1, \dots, a_n \;  \left | 
\begin{array}{ll}
[a_i,a_j] = e   &  2 \leq i, j \leq n, \ \  i \neq j   \\
(a_1 a_k)^2 = (a_k a_1)^2    &  2 \leq k \leq n  \\
(a_1 a_2 \cdots a_n)^2 = e                                              &     
\end{array} 
\right . \Biggr\rangle, 
\end{equation}
where $a_1, \dots, a_n$ are meridians of $Q_1, \dots, Q_n$, respectively.
\end{theorem}

\begin{corollary}\label{bigCn}
The group $\fg{\Cc_n}$ is abelian for $n=1$ and big for $n \geq 2$. 
\end{corollary}

\bigskip

\section{Proof of Theorem \ref{presAn}}\label{pres1}
Take the following affine quadric arrangement $\A_n$, which is composed of $n$ 
quadrics tangent to each other at two points, see Figure \ref{An}. 

\bigskip

\begin{figure}[h]  
\epsfxsize=6cm %width  
\epsfysize=4cm %heigh  
\begin{minipage}{\textwidth}  
\begin{center}  
\epsfbox{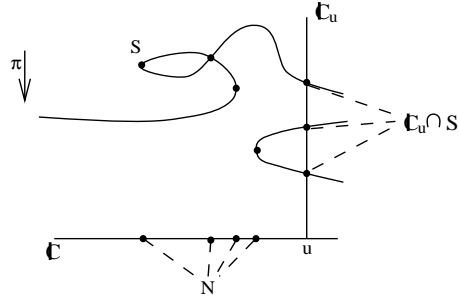}  
\end{center}  
\end{minipage}  
\caption{General setting}\label{setup}  
\end{figure}

In order to find the fundamental group  $\fg{\A_n}$, we need the algorithm of Moishezon-Teicher 
for the global braid monodromy. Consider the following setting (Figure \ref{setup}).  
$S$ is an algebraic curve in $\C^2$,  with $p = \deg(S)$.  Let $\pi: \C^2 \rightarrow \C$ be a 
generic projection  on the first  coordinate. 
Define the fiber $K(x) = \{y \mid (x,y) \in S\}$ in $S$ over a fixed point $x$,  
projected to the $y$-axis. Define $N = \{x \mid \# K(x) < p \}$ and   $M' = \{ s \in S \mid  
\pi_{\mid s} \mbox{ is not \'{e}tale at } s \}$;  note that $\pi (M') = N$.   Let $\{A_j \}^q_{j=1}$ 
be the set of points of $M'$  and $N = \{x_j \}^q_{j=1}$ their projection on the  
$x$-axis.   Recall that $\pi$ is generic, so we assume that $\# (\pi^{-1}(x) \cap M') =1$  
for every $x \in N$. Let $E$ (resp. $D$) be a closed disk on the $x$-axis (resp. the $y$-axis),   
such that  $M' \subset E \times D$ and $N \subset \mbox{Int}(E)$.    
We choose $u \in \partial E$ a real point far enough from the set $N$,  
$x << u$ for every $x \in N$. Define $\C_u = \pi^{-1}(u)$ and number the points of  
$K=\C_u\cap S$ as  $\{1 , \dots , p\}$.  

We construct a g-base for the fundamental group $\pi_1(E - N, u)$.   
Take a set of paths  $\{\gamma_j \}^q_{j=1}$ which connect $u$ with the points  
$\{x_j \}^q_{j=1}$ of $N$.  Now encircle each $x_j$ with a small oriented counterclockwise  
circle $c_j$.  Denote the path segment from $u$ to the boundary of this circle as $\gamma'_j$.   
We define an element (a loop) in the g-base as $\delta_j = {\gamma'}_j c_j {\gamma'}^{-1}_j$.    
Let $B_p[D,K]$ be the braid group, and let $H_1 , \dots , H_{p-1}$ be its frame   
(for complete definitions, see \cite[Section III.2]{16}).    The braid monodromy of $S$  
\cite{2} is a map  $\varphi: \pi_1(E - N, u) \rightarrow B_p[D,K]$ defined as follows:    
every loop in $E - N$ starting at $u$ has liftings to a system of $p$ paths in  
$(E - N) \times D$ starting at each point of $K = 1,\ldots,p$.  Projecting them to $D$  
we get $p$ paths in $D$ defining a motion $\{1(t), \dots , p(t)\}$  
(for $0 \leq t \leq 1$) of $p$ points in $D$ starting and ending at $K$.  This motion  
defines a braid in $B_p [D , K]$.   
By the Artin Theorem \cite{MoTe}, for $j=1, \dots, q$, there exists a halftwist  $Z_j \in  
B_p[D , K]$ and $\epsilon_j \in \Z$, such that $\varphi(\delta_j) = Z_j^{\epsilon_j}$,  
where  $Z_j$ is a halftwist and $\epsilon_j = 1,2$ or $4$ (for an ordinary branch point,  
a node, or a tangency point respectively).  We explain now how to get this $Z_j$.

We explain how to get the braid monodromy around each singularity in $S$.  
Let $A_j$ be a singularity in $S$ and its projection by $\pi$ to the $x$-axis is $x_j$.  
We choose a point $x'_j$ next to $x_j$, such that  $\pi^{-1}(x'_j)$ is a typical fiber.   
We encircle $A_j$ with a  very small circle in a way that the typical  
fiber $\pi^{-1}(x'_j)$ intersects the  circle in two points, say $a ,b$. We fix a skeleton  
$\xi_{x'_j}$ which  connects $a$ and $b$, and denote it as $<a,b>$.   
The Lefschetz diffeomorphism $\Psi$ (see \cite{16}) allows us to  
get a resulting skeleton $(\xi_{x'_j}) \Psi$ in the typical fiber $\C_u$.  
This one defines a motion of its two endpoints. This motion induces a halftwist  
$Z_j = \Delta < (\xi_{x'_j}) \Psi>$. As above,  
$\varphi(\delta_j) = \Delta < (\xi_{x'_j}) \Psi>^{\epsilon_j}$.   
The braid monodromy factorization associated to $S$ is   
$\Delta^2_{p} = \prodl^q_{j=1} \varphi(\delta_j)$.
 
\bigskip
\bigskip

Before proving Theorem \ref{presAn}, we have to study the local monodromy around a common tangency point 
(of $n$ quadrics) and the relation which is derived from this monodromy. This is done in the following propositions.
\begin{proposition}\label{pro14}
Let $O$ be a tangency point, as depicted in Figure \ref{tan1}. 
Let $a_1, a_2, \dots, a_n$ be $n$ points, which are intersections of 
the $n$ quadrics with a typical fiber. 
\begin{figure}[h]
\epsfxsize=11cm %width
\epsfysize=6cm %height
\begin{minipage}{\textwidth}
\begin{center}
\epsfbox{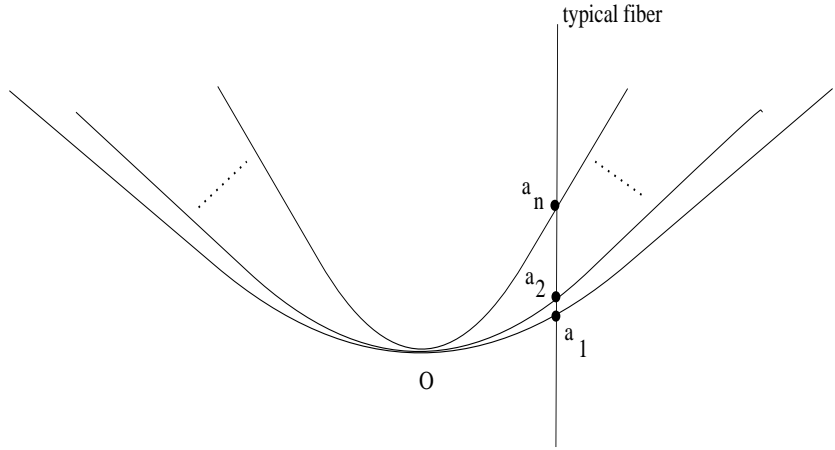}
\end{center}
\end{minipage}
\caption{A tangency point}\label{tan1}
\end{figure}
Then the local monodromy around $O$ is a double fulltwist of  $a_1, a_2, \dots, a_n$ if $O$ is 
one among two tangency points (see e.g. Figure \ref{An}) or a quadruple fulltwist if it is a unique one 
(see e.g. Figure \ref{Bn}).
\end{proposition}

\begin{proof}
Consider $O$ to be locally the tangency point of the quadrics $y=x^2, \dots, y=nx^2$. 
Take a loop $x=e^{2 \pi i t}$ in $y=0$, starting (and ending) 
at some base point and encircling the point $O$, $0 \leq t \leq 1$.  
We have $a_1, a_2, \dots, a_n$ in a typical fiber next to the fiber $F_O$. 

When $t$ is running from $0$ to $1/2$, the point $a_1$ is rotating 
around the other points in a fulltwist $e^{2 \pi i}$, 
the point $a_2$ is rotating along a closed curve which bounds 
a disk, containing the trajectory of the point $a_1$, and so on.  
When $t$ is running from $1/2$ to $1$, we have the same motion. 
This gives us a double fulltwist of the points. 

In a case that $O$ is a unique tangency point, we 
get a doubling of the monodromy according to the B\'ezout Theorem (see e.g. \cite{h}).  
\end{proof}

\begin{proposition}\label{reltan}
If $O$ is one among two tangency points, then we have the relation   
$$
({a_n} \cdots a_2 {a_1})^2 = ({a_1} {a_n} \cdots {a_2})^2 = \dots = ({a_{n-1}} {a_{n-2}} \cdots {a_1} {a_n})^2,  
$$
and if $O$ is a unique tangency point, then we have the relation
$$
({a_n} \cdots a_2 {a_1})^4 = ({a_1} {a_n} \cdots {a_2})^4 = \dots = ({a_{n-1}} {a_{n-2}} \cdots {a_1} {a_n})^4. 
$$
\end{proposition}

\begin{proof}
Say that $O$ is one among two tangency points. The monodromy is a double fulltwist, see Figure \ref{Anelements0}. 
The resulting loops are denoted as  $\tilde a_1, \tilde a_2, \dots, \tilde a_n$. 
\begin{figure}[h]
\epsfxsize=9cm %width
\epsfysize=4cm %height
\begin{minipage}{\textwidth}
\begin{center}
\epsfbox{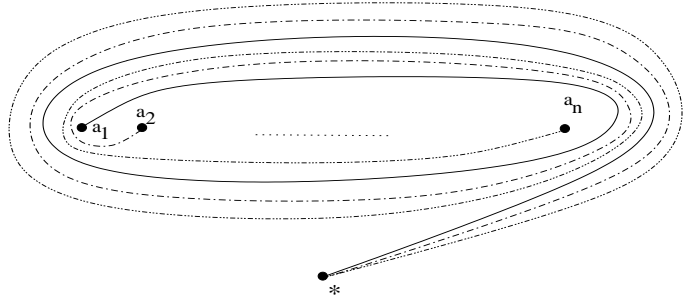}
\end{center}
\end{minipage}
\caption{The resulting loops}\label{Anelements0}
\end{figure}

By the van Kampen Theorem \cite{vankampen}, 
\begin{eqnarray}
{}\tilde{a_1} &=& {a_n} \cdots {a_2} {a_1}  {a_n} \cdots {a_2} {a_1} {a_2}^{-1} \cdots 
{a_{n}}^{-1}  {a_1}^{-1} {a_2}^{-1} \cdots {a_n}^{-1} \nonumber \\ 
{}\tilde{a_2} &=& {a_n} \cdots {a_2} {a_1}  {a_n} \cdots {a_2} {a_1} {a_2} {a_1}^{-1} 
{a_2}^{-1} \cdots {a_{n}}^{-1} {a_1}^{-1}  {a_2}^{-1}  \cdots {a_n}^{-1}\nonumber\\
{}& \dots\dots \nonumber\\
{}\tilde{a_n} &=& {a_n} \cdots {a_1} {a_n} \cdots  {a_1} {a_n} {a_1}^{-1} 
\cdots {a_{n}}^{-1}  {a_1}^{-1} \cdots {a_n}^{-1}  \nonumber,
\end{eqnarray} 
which give us 
\begin{eqnarray*}
{}({a_n} \cdots {a_1})^2  &=& ({a_1} {a_n} \cdots {a_2})^2 \\
{}({a_n} \cdots {a_1})^2  &=& ({a_2} {a_1} {a_n} \cdots {a_3})^2\nonumber \\
{}& \dots\dots \nonumber\\
{}({a_n} \cdots {a_1})^2  &=& ({a_{n-1}} {a_{n-2}} \cdots {a_1} {a_n})^2 \nonumber. 
\end{eqnarray*}

If  $O$ is a unique tangency point, it is clear that we get 
$$
({a_n} \cdots {a_1})^4 = ({a_1} {a_n} \cdots {a_2})^4 = \dots = ({a_{n-1}} {a_{n-2}} \cdots {a_1} {a_n})^4. 
$$
\end{proof}

\begin{figure}[h]  
\epsfxsize=6cm %width  
\epsfysize=5cm %heigh  
\begin{minipage}{\textwidth}  
\begin{center}  
\epsfbox{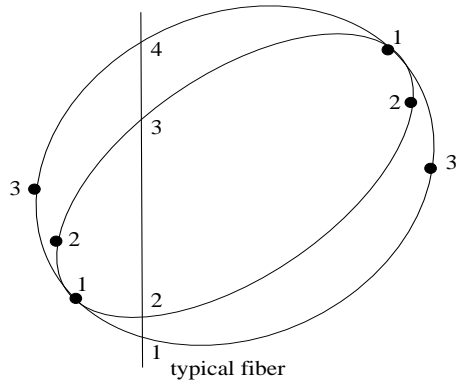}  
\end{center}  
\end{minipage}  
\caption{The arrangement $A_2$}\label{A2fig}  
\end{figure}

\noindent
\underline{\em Proof of Theorem \ref{presAn}:}\\
We are interested in the group $\fg{\A_n}$. We start with $\A_1$ (a smooth quadric), 
and it is easy to see that the group $\fg{\A_1}$ is $\Z_2$. Let us consider the arrangement $\A_2$, see Figure \ref{A2fig}.
We want to compute the braid monodromy factorization of  $\A_2$. Consider the above general setting 
(see Figure \ref{setup}).
Let $K=\{1, 2, 3, 4\}$ and let $\{j\}^{3}_{j=1}$ be singular points of $\pi_1$ on the right side of 
the chosen typical fiber as follows:   $1$ is a tangency point,  $2$ and $3$ are  branch points of the quadrics.

We are looking for $\varphi(\delta_j)$ for $j = 1, 2, 3$.   
So we choose a $g$-base $\{\delta_j\}^{3}_{j=1}$ of $\pi_1(E - N,u)$, such that  
each $\delta_j$ is constructed from a path $\gamma_j$ below the real line 
and a counterclockwise small circle around the points in $N$. 
 
The diffeomorphisms which are induced from passing through branch points were defined in  
\cite{MoTe} as $\Delta^{\frac{1}{2}}_{I_2\R}<k>$ and  $\Delta^{\frac{1}{2}}_{I_4 I_2}<k>$. 
We recall the precise  definition from \cite{MoTe}.  
Consider a typical fiber on the left side of a branch point 
(locally defined by $y^2-x=0$). The typical fiber intersects the quadric in two 
complex points. Passing through this point, the two complex points move to the k'th place and rotate  
in a counterclockwise $90^o$ twist. They become real and numbered as $k, k+1$.

Find first the skeleton $\xi_{x'_j}$ related to each singular point. 
Then compute the local diffeomorphisms $\delta_{j}$  induced from singular points $j$.  
The monodromy table is:
\begin{center} 
\begin{tabular}{cccc} \\ 
$j$ & $\xi_{x'_j}$ & $\epsilon_{j} 
$ & $\delta_{j}$ \\ \hline 
1 & $<3,4>$ & 4 & $\Delta^2<3,4>$\\ 
2 & $<2,3>$ & 1 & $\Delta^{\frac{1}{2}}_{I_2\R}<2>$\\ 
3 & $<1,2>$ & 1 & $\Delta^{\frac{1}{2}}_{I_4 I_2}<1>$ 
\end{tabular} 
\end{center} 
Using \cite{MoTe},  we compute the skeleton  $(\xi_{x'_j}) \Psi_{\gamma'_j}$ to each $j$ by applying 
to the  skeleton $\xi_{x'_j}$ the product $\prodl_{i=j-1}^{1} \delta_{i}$.  
\begin{list}{}{\leftmargin 2cm \labelsep .4cm \labelwidth 3cm} 
\item $(\xi_{x'_1}) \Psi_{\gamma'_1} =<3,4>= z_{3\ 4}$\\ 
$\varphi(\delta_1) = Z^4_{3\ 4}$ 
\vspace{-.3cm} 
\begin{figure}[h!] 
\epsfysize=0.4cm %height 
\begin{minipage}{\textwidth} 
\begin{center} 
\epsfbox{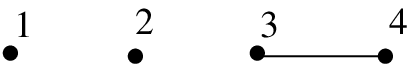} 
\end{center} 
\end{minipage} 
\end{figure} 
\item $(\xi_{x'_2}) \Psi_{\gamma'_2} = <2,3> \Delta^2<3,4>= z_{2\ 3 }^{Z^2_{3\ 4}}$\\ 
$\varphi(\delta_2) = {Z_{2\ 3 }}^{Z^2_{3\ 4}}$ 
\begin{figure}[h!] 
\epsfysize=1.1cm %height 
\begin{minipage}{\textwidth} 
\begin{center} 
\epsfbox{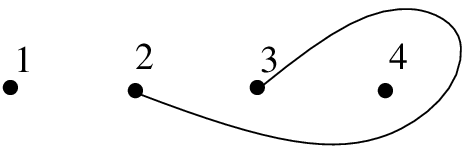} 
\end{center} 
\end{minipage} 
\end{figure} 
\item $(\xi_{x'_3}) \Psi_{\gamma'_3}=<1,2> \Delta^{\frac{1}{2}}_{I_2 \R}<2>\Delta^2<3,4>={z}_{1\ 4}$\\ 
$\varphi(\delta_3) = {Z}_{1\ 4}$ 
\begin{figure}[h!] 
\epsfysize=0.8cm %height 
\begin{minipage}{\textwidth} 
\begin{center} 
\epsfbox{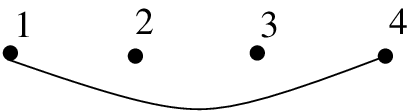} 
\end{center} 
\end{minipage} 
\end{figure} 
\end{list} 
On the left side of the typical fiber (Figure \ref{A2fig}) we have three singularities and 
the related monodromy table is: 
\begin{center} 
\begin{tabular}{cccc} \\ 
$j$ & $\xi_{x'_j}$ & $\epsilon_{j} 
$ & $\delta_{j}$ \\ \hline 
1 & $<1,2>$ & 4 & $\Delta^2<1,2>$\\ 
2 & $<2,3>$ & 1 & $\Delta^{\frac{1}{2}}_{I_2\R}<2>$\\ 
3 & $<1,2>$ & 1 & $\Delta^{\frac{1}{2}}_{I_4 I_2}<1>$ 
\end{tabular} 
\end{center} 
And again by the Moishezon-Teicher  algorithm \cite{MoTe} we have the following resulting braids:
\begin{list}{}{\leftmargin 2cm \labelsep .4cm \labelwidth 3cm} 
\item $(\xi_{x'_1}) \Psi_{\gamma'_1} =<1,2>= z_{1\ 2}$\\ 
$\varphi(\delta_1) = Z^4_{1\ 2}$ 
\vspace{-.3cm} 
\begin{figure}[h!] 
\epsfysize=0.4cm %height 
\begin{minipage}{\textwidth} 
\begin{center} 
\epsfbox{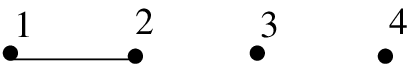} 
\end{center} 
\end{minipage} 
\end{figure} 
\item $(\xi_{x'_2}) \Psi_{\gamma'_2} = <2,3> \Delta^2<1,2>= z_{2\ 3 }^{Z^2_{1\ 2}}$\\ 
$\varphi(\delta_2) = {Z_{2\ 3 }}^{Z^2_{1\ 2}}$ 
\begin{figure}[h!] 
\epsfysize=1.2cm %height 
\begin{minipage}{\textwidth} 
\begin{center} 
\epsfbox{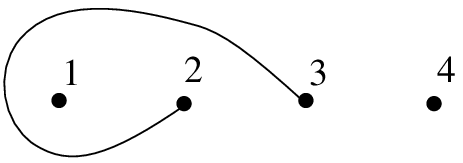} 
\end{center} 
\end{minipage} 
\end{figure} 
\item $(\xi_{x'_3}) \Psi_{\gamma'_3}=<1,2> \Delta^{\frac{1}{2}}_{I_2 \R}<2>\Delta^2<1,2>=\bar{z}_{1\ 4}$\\ 
$\varphi(\delta_3) = \bar{Z}_{1\ 4}$ 
\begin{figure}[h!]
\epsfysize=0.7cm %height 
\begin{minipage}{\textwidth} 
\begin{center} 
\epsfbox{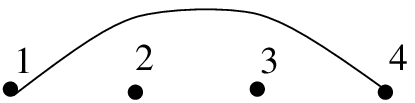} 
\end{center} 
\end{minipage} 
\end{figure} 
\end{list} 

According to the van Kampen Theorem \cite{vankampen}, we get the following set of relations for  
$\fg{\A_2}$: 
\begin{eqnarray*}
{}&&(3 4)^2 = (4 3)^2\\
{}&&2 = 4 3 4^{-1}\\
{}&&1 = 4\\
{}&&(1 2)^2 = (2 1)^2\\
{}&&3 = 2 1 2 1^{-1} 2^{-1}\\
{}&&4 = 3 2 1  2^{-1} 3^{-1}\\
{}&&4321 = e \ \ \mbox{(the projective relation).}
\end{eqnarray*}
By an easy simplification, we get  $\fg{\A_2} \simeq \langle a_1, a_2   \left | 
 (a_1 a_2)^2 = e \right \rangle$.

Now we consider the arrangement $\A_3$. In a similar way as done for $\A_2$, we compute 
the relevant monodromy and we are able to find out that $\fg{\A_3} \simeq \langle a_1, a_2, a_3 
  \left |  (a_1 a_2 a_3)^2 = e \right \rangle$. 

In order to compute the group $\fg{\A_n}$, we have to apply again the Moishezon-Teicher algorithm 
for the global monodromy. On each one of the sides of the typical fiber in Figure \ref{An} we have 
one tangency point and $n$ branch points. 
The relation related to a tangency point was proved in Proposition \ref{reltan}. Since we have two 
tangency points in the arrangement, we have relations (\ref{an1}) and (\ref{an5}).  
The relations relating to the right branch points are (\ref{an2})- -(\ref{an4}), and 
the relations   relating to the left branch points are (\ref{an6})- -(\ref{an8}): 
\begin{eqnarray}
\ && (a_{n+1} a_{n+2} \cdots a_{2n})^2 =(a_{2n} a_{n+1} \cdots a_{2n-1})^2 = \dots 
= (a_{n+2} a_{n+3} \cdots a_{2n} a_{n+1})^2 \label{an1}\\
\  && a_{n}= a_{2n} a_{2n-1} \cdots a_{n+2} a_{n+1} a_{n+2}^{-1} \cdots a_{2n-1}^{-1} a_{2n}^{-1}  
\label{an2}\\
\  && a_{n-1}= a_{2n} a_{2n-1} \cdots a_{n+3} a_{n+2} a_{n+3}^{-1} \cdots a_{2n-1}^{-1} a_{2n}^{-1}  
\label{an3}\\
\  && a_{n-2}= a_{2n} a_{2n-1} \cdots a_{n+4} a_{n+3} a_{n+4}^{-1} \cdots a_{2n-1}^{-1} a_{2n}^{-1}  
\label{an3a}\\
\  &&\dots\dots \nonumber \\
\  && a_{1}= a_{2n} \label{an4}\\
\  && (a_{1} a_{2} \cdots a_{n})^2 =(a_{n} a_{1} \cdots a_{n-1})^2 = \dots= 
(a_{2} a_{3} \cdots a_{n} a_{1})^2  \label{an5}\\
\  && a_{n+1}= a_{n} a_{n-1} \cdots a_{1} a_{n} a_{1}^{-1} \cdots a_{n-1}^{-1} a_{n}^{-1} \label{an6} \\
\  && a_{n+2}= a_{n+1} a_{n} \cdots a_{1} a_{n-1} a_{1}^{-1} \cdots a_{n}^{-1} a_{n+1}^{-1} \label{an7} \\
\  && a_{n+3}= a_{n+2} a_{n+1} \cdots a_{1} a_{n-2} a_{1}^{-1} \cdots a_{n+1}^{-1} a_{n+2}^{-1} 
\label{an7a} \\
\  && \dots\dots \nonumber\\
\  && a_{2n}= a_{2n-1} a_{2n-2} \cdots a_{2} a_{1} a_{2}^{-1} \cdots a_{2n-2}^{-1} a_{2n-1}^{-1} 
\label{an8} \\
\  && a_{2n} a_{2n-1} a_{2n-2} \cdots a_{2} a_{1} = e \  \ \mbox{(the projective relation).} \label{an9}.
\end{eqnarray}
We substitute relations  (\ref{an2})- -(\ref{an4}) in  relations (\ref{an6})- -(\ref{an8})  
and get for each one of them one of the following equations:
$$
(a_{n+1} a_{n+2} \cdots a_{2n})^{-2} =(a_{2n} a_{n+1} \cdots a_{2n-1})^{-2} = \dots 
= (a_{n+2} a_{n+3} \cdots a_{2n} a_{n+1})^{-2}.
$$
These equations appear already in (\ref{an1}). Therefore  relations (\ref{an6})- -(\ref{an8}) 
are redundant. 

Also relation (\ref{an5}) gets the same form of (\ref{an1}), again by substituting 
relations  (\ref{an2})- -(\ref{an4}), and therefore it is redundant too.

By the same substitution, relation  (\ref{an9}) is rewritten as $(a_{n+1} a_{n+2} \cdots a_{2n})^2 =~e$. 
Hence, we get the presentation (\ref{equAn}), as needed.
\hfill $\Box$

\bigskip

\noindent
\underline{\em Proof of Corollary \ref{bigAn}:}\\
We computed above the group  $\fg{\A_2}$. This group is isomorphic to $\Z * \Z_2$, 
which is known to be big. Since the group $\fg{\A_2}$ is a quotient of  $\fg{\A_n}$  
for $n \geq 2$,  the groups $\fg{\A_n}$  are big too.
\hfill $\Box$

\bigskip

\begin{remark}
There is another way to prove that the groups $\fg{\A_n}$ are big (for $n \geq 2$). 
These groups are isomorphic to $\Z^{n-1} * \Z_2$, which contain a non-abelian free subgroup.
\end{remark}

\bigskip

\section{Proof of Theorem \ref{presBn}}\label{pres2}
\begin{figure}[h]  
\epsfxsize=6cm %width  
\epsfysize=5cm %heigh  
\begin{minipage}{\textwidth}  
\begin{center}  
\epsfbox{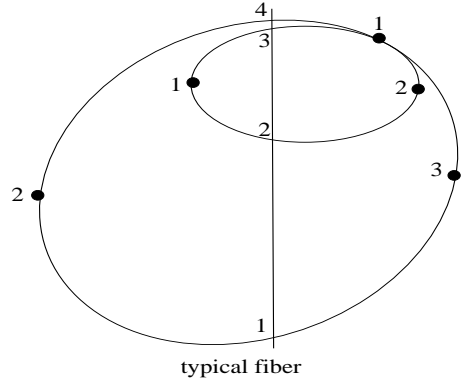}  
\end{center}  
\end{minipage}  
\caption{The arrangement $B_2$}\label{B2fig}  
\end{figure}

The arrangement $\B_n$ is depicted in Figure \ref{Bn}. 
We start with $\B_1$ (a smooth quadric), and it is easy to see that the group $\fg{\B_1}$ is $\Z_2$. 
Now we compute the braid monodromy factorization of the curve $\B_2$. We follow Figure \ref{B2fig}.
Let $\pi_1 : E \times D \rightarrow E$ be the projection to $E$. 
Take $u \in \partial E$, such that  $\C_u$ is a typical fiber and $M \in \C_u$ is a real point. 
Let $K = K(M)=\{1, 2, 3, 4\}$. 
Let $\{j\}^{3}_{j=1}$ be singular points of $\pi_1$ on the right side of the chosen typical fiber 
as follows:   $1$ is a tangency point,  $2$ and $3$ are  branch points of the quadrics. 
Let $N = \{x(j) = x_j \mid 1 \le j \le 3 \}$, such that $N \subset E - \partial E, N \subset E$. 
We start with the monodromy table:
\begin{center} 
\begin{tabular}{cccc} \\ 
$j$ & $\xi_{x'_j}$ & $\epsilon_{j} 
$ & $\delta_{j}$ \\ \hline 
1 & $<3,4>$ & 8 & $\Delta^4<3,4>$\\ 
2 & $<2,3>$ & 1 & $\Delta^{\frac{1}{2}}_{I_2\R}<2>$\\ 
3 & $<1,2>$ & 1 & $\Delta^{\frac{1}{2}}_{I_4 I_2}<1>$ 
\end{tabular} 
\end{center} 
By the algorithm from  \cite{MoTe},  we compute the  monodromies: 
\begin{list}{}{\leftmargin 2cm \labelsep .4cm \labelwidth 3cm} 
\item $(\xi_{x'_1}) \Psi_{\gamma'_1} =<3,4>= z_{3\ 4}$\\ 
$\varphi(\delta_1) = Z^8_{3\ 4}$ 
\vspace{-.3cm} 
\begin{figure}[h!]
\epsfysize=0.4cm %height 
\begin{minipage}{\textwidth} 
\begin{center} 
\epsfbox{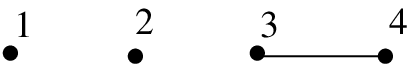} 
\end{center} 
\end{minipage} 
\end{figure} 
\item $(\xi_{x'_2}) \Psi_{\gamma'_2} = <2,3> \Delta^4<3,4>= z_{2\ 3 }^{Z^4_{3\ 4}}$\\ 
$\varphi(\delta_2) = {Z_{2\ 3 }}^{Z^4_{3\ 4}}$ 
\begin{figure}[h!]
\epsfysize=1.6cm %height 
\begin{minipage}{\textwidth} 
\begin{center} 
\epsfbox{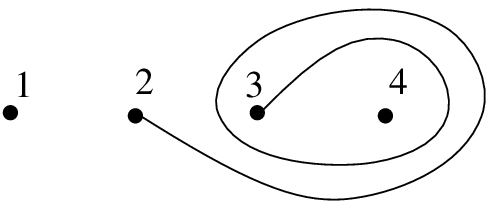} 
\end{center} 
\end{minipage} 
\end{figure} 
\item $(\xi_{x'_3}) \Psi_{\gamma'_3}=<1,2> \Delta^{\frac{1}{2}}_{I_2 \R}<2>\Delta^4<3,4>=
{z}_{1\ 4}^{Z^2_{3\ 4}}$\\ 
$\varphi(\delta_3) = {{Z}_{1\ 4}}^{Z^2_{3\ 4}}$ 
\begin{figure}[h!]
\epsfysize=1.5cm %height 
\begin{minipage}{\textwidth} 
\begin{center} 
\epsfbox{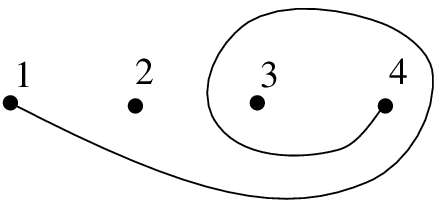} 
\end{center} 
\end{minipage} 
\end{figure} 
\end{list} 
On the left side of the typical fiber in Figure \ref{B2fig} we have two branch points and 
the related monodromy table is: 
\begin{center} 
\begin{tabular}{cccc} \\ 
$j$ & $\xi_{x'_j}$ & $\epsilon_{j} 
$ & $\delta_{j}$ \\ \hline
1 & $<2,3>$ & 1 & $\Delta^{\frac{1}{2}}_{I_2\R}<2>$\\
2 & $<1,2>$ & 1 & $\Delta^{\frac{1}{2}}_{I_4 I_2}<1>$ 
\end{tabular} 
\end{center} 
And we have the following resulting braids:
\begin{list}{}{\leftmargin 2cm \labelsep .4cm \labelwidth 3cm} 
\item $(\xi_{x'_1}) \Psi_{\gamma'_1} =<2,3>= z_{2\ 3}$\\ 
$\varphi(\delta_1) = Z_{2\ 3}$ 
\vspace{-.3cm} 
\begin{figure}[h!]
\epsfysize=0.4cm %height 
\begin{minipage}{\textwidth} 
\begin{center} 
\epsfbox{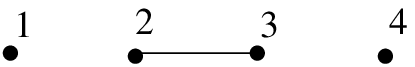} 
\end{center} 
\end{minipage} 
\end{figure} 
\item $(\xi_{x'_2}) \Psi_{\gamma'_2} = <1,2> \Delta^{\frac{1}{2}}_{I_2\R}<2>= z_{1\ 4}^{Z^{2}_{1\ 3}}$\\ 
$\varphi(\delta_2) = {Z_{1\ 4}}^{Z^{2}_{1\ 3}}$ 
\begin{figure}[h!]
\epsfysize=1.2cm %height 
\begin{minipage}{\textwidth} 
\begin{center} 
\epsfbox{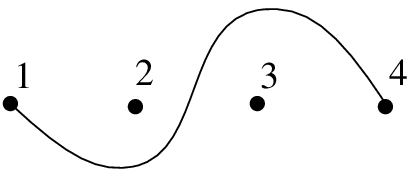} 
\end{center} 
\end{minipage} 
\end{figure} 
\end{list} 
The relations which are related to these braids are as follows:
\begin{eqnarray*}
{}&&(3 4)^4 = (4 3)^4\\
{}&&2 = 4 3 4 3 4^{-1} 3^{-1} 4^{-1} \\
{}&&1 = 4 3 4 3^{-1} 4^{-1} \\
{}&&2 = 3\\
{}&&1 = 3^{-1} 4 3\\
{}&&4321 = e \ \ \mbox{(the projective relation)}.
\end{eqnarray*}
By an easy simplification, we get  $\fg{\B_2} \simeq \langle a_1, a_2   \left | 
 (a_1 a_2)^2 = e \right \rangle$.

In a similar way, we obtain the group $\fg{\B_3} \simeq \langle a_1, a_2, a_3 
  \left |  (a_1 a_2 a_3)^2 = e \right \rangle$. 

In order to compute the group $\fg{\B_n}$, we have to apply again the Moishezon-Teicher algorithm 
for the global monodromy. On the right side of the typical fiber in Figure \ref{Bn} we have 
one tangency point and $n$ branch points. On the left side we have  $n$ branch points. 
The relations related to the right points are:
\begin{eqnarray}
\ &&  (a_{n+1} a_{n+2} \cdots a_{2n})^4 = (a_{2n} a_{n+1} \cdots  a_{2n-1})^4 = \dots = 
(a_{n+2} a_{n+3} \cdots a_{2n} a_{n+1})^4  \label{bn1}\\
\  && a_{n}= a_{2n} a_{2n-1} \cdots a_{n+2} a_{n+1} a_{2n} a_{2n-1} \cdots 
a_{n+2} a_{n+1}  a_{n+2}^{-1} \cdots  a_{2n-1}^{-1} a_{2n}^{-1} \cdot \label{bn2}\\
\  &&\hspace{1cm}  \cdot a_{n+1}^{-1} a_{n+2}^{-1} \cdots a_{2n-1}^{-1} a_{2n}^{-1} 
\nonumber \\
\  && a_{n-1}= a_{2n} a_{2n-1} \cdots a_{n+2}  a_{n+1} a_{2n} a_{2n-1} \cdots a_{n+3} a_{n+2} 
a_{n+3}^{-1} \cdots a_{2n-1}^{-1} a_{2n}^{-1}  \cdot \label{bn3}\\
\  &&\hspace{1cm}  \cdot a_{n+1}^{-1} a_{n+2}^{-1} \cdots  a_{2n-1}^{-1}a_{2n}^{-1} \nonumber \\
\  && \dots\dots  \nonumber \\
\  && a_{1}= a_{2n} a_{2n-1} \cdots a_{n+2} a_{n+1} a_{2n} a_{n+1}^{-1} a_{n+2}^{-1} \cdots  
a_{2n-1}^{-1} a_{2n}^{-1} \label{bn4},
\end{eqnarray}
and the relations related to the left points are:
\begin{eqnarray}
\  && a_{n}= a_{n+1} \label{bn5}\\
\  &&a_{n-1}= a_{n+1}^{-1}  a_{n+2} a_{n+1} \label{bn6}\\
\   &&\dots\dots  \nonumber \\
\   && a_{1}= a_{n+1}^{-1} a_{n+2}^{-1} \cdots  a_{2n-1}^{-1} a_{2n} a_{2n-1} \cdots a_{n+2} 
a_{n+1}\label{bn7}.
\end{eqnarray}
The projective relation is:
\begin{eqnarray}
\  a_{2n} a_{2n-1} a_{2n-2} \cdots a_{2} a_{1} = e. \label{bn8}
\end{eqnarray}
We substitute relation (\ref{bn5}) in relation (\ref{bn2}) and get 
$$(a_{n+1} a_{2n} a_{2n-1} \cdots a_{n+2})^2 = (a_{2n} a_{2n-1} \cdots a_{n+1})^2.$$
Substituting relation (\ref{bn6}) in relation (\ref{bn3}), we get  
$$(a_{n+1} a_{2n} a_{2n-1} \cdots a_{n+2})^2 = (a_{n+2} a_{n+1}  a_{2n} \cdots a_{n+3})^2.$$
Repeating this procedure, we get for each one of the substitutions (the last step is 
the substitution of (\ref{bn7})  in (\ref{bn4})) one of the following equations:  
$$
(a_{2n} a_{2n-1} \cdots a_{n+1})^{2} =(a_{n+1} a_{2n} \cdots a_{n+2})^{2} = \dots 
= (a_{2n-1} \cdots a_{n+1} a_{2n})^{2}.
$$
Therefore relation (\ref{bn1}) is redundant. 

The projective relation is rewritten as $(a_{n+1} a_{n+2} \cdots a_{2n})^2 =e$ (by the same 
substitutions). Hence we get the presentation (\ref{equBn}).
\hfill $\Box$

\bigskip

\noindent
\underline{\em Proof of Corollary \ref{bigBn}:}\\
The proof is the same one  as in Corollary \ref{bigAn}.
\hfill $\Box$

\bigskip

\section{Proof of Theorem \ref{presCn}}\label{pres3}
Let $\Cc_n$ be a quadric arrangement, composed of $n$ quadrics as shown in Figure 
\ref{Cn}. Each one of the quadrics $Q_2, \dots, Q_n$ is tangent to the quadric 
$Q_1$ at one tangency point and intersects each one of the other quadrics at four points.

We start with the arrangement $\Cc_1$. The group $\fg{\Cc_1}$ is again $\Z_2$. 

\begin{figure}[h]  
\epsfxsize=6cm %width  
\epsfysize=4cm %heigh  
\begin{minipage}{\textwidth}  
\begin{center}  
\epsfbox {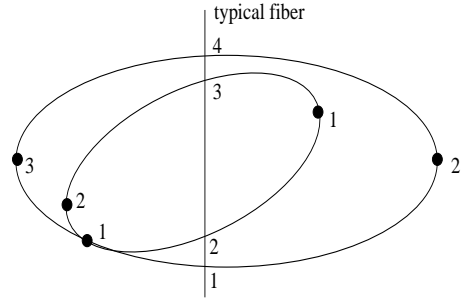}  
\end{center}  
\end{minipage}  
\caption{The arrangement $C_2$}\label{C2fig}  
\end{figure}

Let us consider now the arrangement $\Cc_2$ from Figure \ref{C2fig}.  
We want to compute the braid monodromy factorization of $\Cc_2$.  
Let $\pi_1 : E \times D \rightarrow E$ be the projection to $E$. 
Take $u \in \partial E$, such that  $\C_u$ is a typical fiber and $M \in \C_u$ is a real point. 
Let $K = K(M)=\{1, 2, 3, 4\}$. 
Let $\{j\}^{2}_{j=1}$ be two branch points of $\pi_1$ on the right side of the chosen typical fiber. 
Their monodromy table and braids are:
\begin{center} 
\begin{tabular}{cccc} \\ 
$j$ & $\xi_{x'_j}$ & $\epsilon_{j} 
$ & $\delta_{j}$ \\ \hline 
1 & $<2,3>$ & 1 & $\Delta^{\frac{1}{2}}_{I_2\R}<2>$\\ 
2 & $<1,2>$ & 1 & $\Delta^{\frac{1}{2}}_{I_4 I_2}<1>$ 
\end{tabular} 
\end{center} 
\begin{list}{}{\leftmargin 2cm \labelsep .4cm \labelwidth 3cm} 
\item $(\xi_{x'_1}) \Psi_{\gamma'_1} =<2,3>= z_{2\ 3}$\\ 
$\varphi(\delta_1) = Z_{2\ 3}$ 
\vspace{-.3cm} 
\begin{figure}[h!]
\epsfysize=0.4cm %height 
\begin{minipage}{\textwidth} 
\begin{center} 
\epsfbox{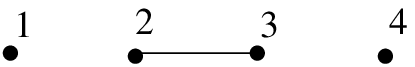} 
\end{center} 
\end{minipage} 
\end{figure} 
\item $(\xi_{x'_2}) \Psi_{\gamma'_2} = <1,2> \Delta^{\frac{1}{2}}_{I_2\R}<2> = 
z_{1\ 4}^{Z^{2}_{1\ 3}}$\\ 
$\varphi(\delta_2) = {Z_{1\ 4}}^{Z^{2}_{1\ 3}}$ 
\begin{figure}[h!]
\epsfysize=1.3cm %height 
\begin{minipage}{\textwidth} 
\begin{center} 
\epsfbox{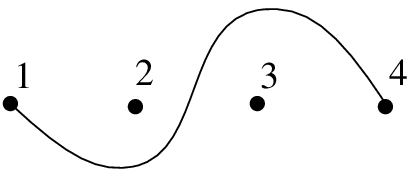} 
\end{center} 
\end{minipage} 
\end{figure} 
\end{list} 
On the left side of the typical fiber we have $1$ which is a tangency point, and 
$2$, $3$ which are branch points. According to Proposition \ref{reltan}, 
since $1$ is the unique tangency point, the monodromy related to it is a quadruple fulltwist. 
The monodromy table and resulting braids are:
\begin{center} 
\begin{tabular}{cccc} \\ 
$j$ & $\xi_{x'_j}$ & $\epsilon_{j} 
$ & $\delta_{j}$ \\ \hline 
1 & $<1,2>$ & 8 & $\Delta^4<1,2>$\\ 
2 & $<2,3>$ & 1 & $\Delta^{\frac{1}{2}}_{I_2\R}<2>$\\ 
3 & $<1,2>$ & 1 & $\Delta^{\frac{1}{2}}_{I_4 I_2}<1>$ 
\end{tabular} 
\end{center} 
\begin{list}{}{\leftmargin 2cm \labelsep .4cm \labelwidth 3cm} 
\item $(\xi_{x'_1}) \Psi_{\gamma'_1} =<1,2>= z_{1\ 2}$\\ 
$\varphi(\delta_1) = Z^8_{1\ 2}$ 
\vspace{-.3cm} 
\begin{figure}[h!] 
\epsfysize=0.4cm %height 
\begin{minipage}{\textwidth} 
\begin{center} 
\epsfbox{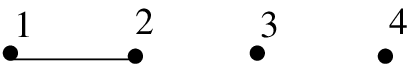} 
\end{center} 
\end{minipage} 
\end{figure} 
\item $(\xi_{x'_2}) \Psi_{\gamma'_2} = <2,3> \Delta^4<1,2>= z_{2\ 3 }^{Z^4_{1\ 2}}$\\ 
$\varphi(\delta_2) = {Z_{2\ 3 }}^{Z^4_{1\ 2}}$ 
\begin{figure}[h!] 
\epsfysize=1.5cm %height 
\begin{minipage}{\textwidth} 
\begin{center} 
\epsfbox{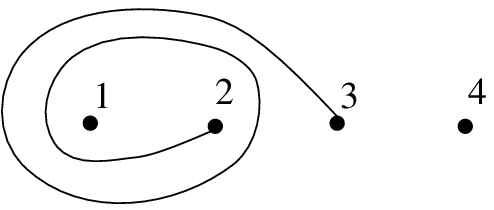} 
\end{center} 
\end{minipage} 
\end{figure} 
\item $(\xi_{x'_3}) \Psi_{\gamma'_3}=<1,2> \Delta^{\frac{1}{2}}_{I_2 \R}<2>\Delta^4<1,2>=\bar{z}_{1\ 4}^{Z^2_{1\ 2}}$\\ 
$\varphi(\delta_3) = {{{\bar{Z}}_{1\ 4}}}^{Z^2_{1\ 2}}$ 
\begin{figure}[h!]
\epsfysize=1.3cm %height 
\begin{minipage}{\textwidth} 
\begin{center} 
\epsfbox{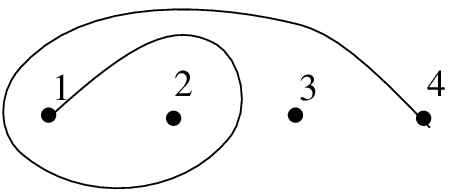} 
\end{center} 
\end{minipage} 
\end{figure} 
\end{list} 
The list of relations for the group $\fg{\Cc_2}$ is
\begin{eqnarray*}
{}&&2  = 3\\
{}&&1 = 3^{-1} 4 3\\
{}&&(1 2)^4  = (2 1)^4\\
{}&&3 = 2 1 2 1 2 1^{-1} 2^{-1} 1^{-1}  2^{-1} \\
{}&&3^{-1} 4 3 = 2 1 2 1 2^{-1} 1^{-1}  2^{-1} \\
{}&&4321 = e \ \ \mbox{(the projective relation)}.
\end{eqnarray*}
By simplification, we get  $\fg{\Cc_2} \simeq \langle a_1, a_2   \left | (a_1 a_2)^2 = e \right \rangle$.

\begin{figure}[h]  
\epsfxsize=7cm %width  
\epsfysize=5cm %heigh  
\begin{minipage}{\textwidth}  
\begin{center}  
\epsfbox {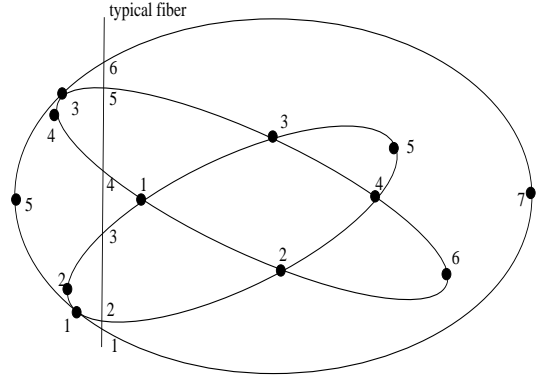}  
\end{center}  
\end{minipage}  
\caption{The arrangement $C_3$}\label{C3fig}  
\end{figure}

The arrangement $\Cc_3$ is depicted in Figure \ref{C3fig}.
Here $K = K(M)=\{1, 2, 3, 4, 5, 6\}$. 
On the right side of the chosen typical fiber there are four  nodes and three branch points. 
The monodromy table and braids are:
\begin{center} 
\begin{tabular}{cccc} \\ 
$j$ & $\xi_{x'_j}$ & $\epsilon_{j} 
$ & $\delta_{j}$ \\ \hline 
1 & $<3,4>$ & 2 & $\Delta<3,4>$\\ 
2 & $<2,3>$ & 2 & $\Delta<2,3>$\\ 
3 & $<4,5>$ & 2 & $\Delta<4,5>$\\ 
4 & $<3,4>$ & 2 & $\Delta<3,4>$\\ 
5 & $<4,5>$ & 1 & $\Delta^{\frac{1}{2}}_{I_2 \R}<4>$\\ 
6 & $<2,3>$ & 1 & $\Delta^{\frac{1}{2}}_{I_4 I_2}<2>$\\
7 & $<1,2>$ & 1 & $\Delta^{\frac{1}{2}}_{I_6 I_4}<1>$
\end{tabular} 
\end{center} 
\begin{list}{}{\leftmargin 2cm \labelsep .4cm \labelwidth 3cm} 
\item $\varphi(\delta_1) = Z^2_{3\ 4}$ \ \ \ \ \ \  $\vcenter{\hbox{\epsfbox{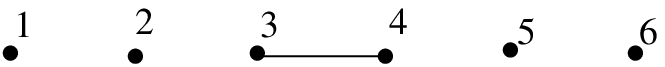}}}$
\item $\varphi(\delta_2) = {Z^2_{2\ 4}}$ \ \ \ \ \ \  $\vcenter{\hbox{\epsfbox{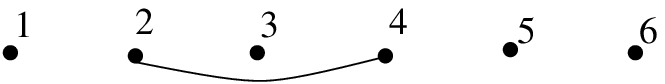}}}$
\item $\varphi(\delta_3) = {{Z}^2_{3\ 5}}^{Z^2_{3\ 4}}$  \ $\vcenter{\hbox{\epsfbox{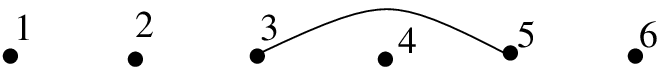}}}$
\item $\varphi(\delta_4) = {{Z}^2_{2\ 5}}^{Z^2_{2\ 4}}$ \ $\vcenter{\hbox{\epsfbox{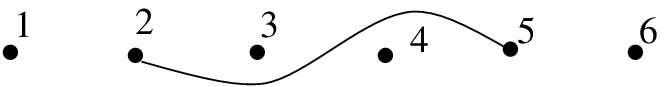}}}$
\item $\varphi(\delta_5) = {Z_{2\ 3}}$ \ \ \ \ \ $\vcenter{\hbox{\epsfbox{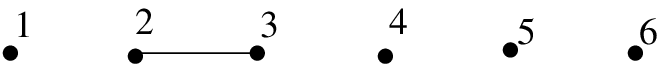}}}$
\item $\varphi(\delta_6) = {Z_{4\ 5}}$  \ \ \ \ \  $\vcenter{\hbox{\epsfbox{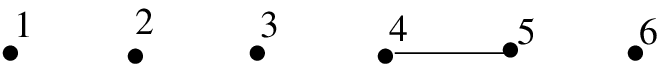}}}$
\item $\varphi(\delta_7) = {\tilde{Z}_{1\ 6}}$ \ \ \  \  \  $\vcenter{\hbox{\epsfbox{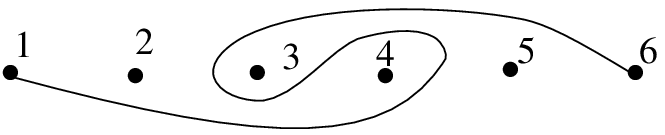}}}$
\end{list} 
On the left side of the typical fiber we have five more singularities:  two of them are 
tangency points and the rest are branch points.
\begin{center} 
\begin{tabular}{cccc} \\ 
$j$ & $\xi_{x'_j}$ & $\epsilon_{j} 
$ & $\delta_{j}$ \\ \hline 
1 & $<1,2>$ & 8 & $\Delta^4<1,2>$\\ 
2 & $<2,3>$ & 1 & $\Delta^{\frac{1}{2}}_{I_2 \R}<2>$\\ 
3 & $<3,4>$ & 8 & $\Delta^4<3,4>$\\ 
4 & $<2,3>$ & 1 & $\Delta^{\frac{1}{2}}_{I_4 I_2}<2>$\\
5 & $<1,2>$ & 1 & $\Delta^{\frac{1}{2}}_{I_6 I_4}<1>$
\end{tabular} 
\end{center} 
\begin{list}{}{\leftmargin 2cm \labelsep .4cm \labelwidth 3cm} 
\item $\varphi(\delta_1) = Z^8_{1\ 2}$ \ \ \ \ \ \ \ \ $\vcenter{\hbox{\epsfbox{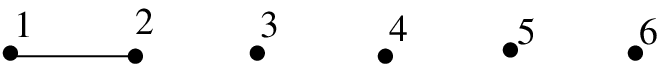}}}$
\item $\varphi(\delta_2) = {Z_{2\ 3}}^{Z^4_{1\ 2}}$   $\vcenter{\hbox{\epsfbox{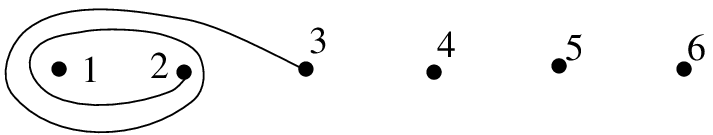}}}$
\item $\varphi(\delta_3) = {Z}^8_{5\ 6}$  \ \ \ \ \ \ \ $\vcenter{\hbox{\epsfbox{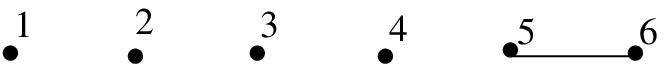}}}$
\item $\varphi(\delta_4) = {{Z}_{4\ 5}}^{Z^4_{5\ 6}}$   \ \ \
$\vcenter{\hbox{\epsfbox{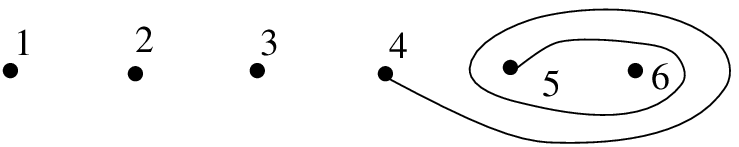}}}$
\item $\varphi(\delta_5) = {\tilde{Z}_{1\ 6}}$  \ \ \ \ \ \ 
$\vcenter{\hbox{\epsfbox{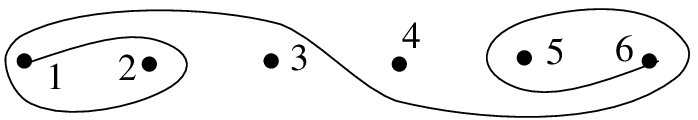}}}$
\end{list} 
The group  $\fg{\Cc_3}$ is then generated by $1, \dots, 6$ and admits the following presentation:
\begin{eqnarray*}
{}&&[3,4]=e\\
{}&&[2,4]=e\\
{}&&[4^{-1} 3 4, 5]=e\\
{}&&[2, 4^{-1} 5 4]=e\\
{}&& 2  = 3\\
{}&& 4  = 5\\
{}&&6  = 5 4 3 4^{-1} 1 4 3^{-1} 4^{-1} 5^{-1}\\
{}&&(1 2)^4 = (2 1)^4\\
{}&&3 = 2 1 2 1 2 1^{-1} 2^{-1} 1^{-1} 2^{-1}\\
{}&&(5 6)^4=(6 5)^4\\
{}&&4 = 6 5 6 5 6^{-1} 5^{-1} 6^{-1}\\
{}&&3 2 1 2 1 2^{-1} 1^{-1} 2^{-1} 3^{-1} = 6 5 6  5^{-1} 6^{-1}\\
{}&&654321 = e \ \ \mbox{(the projective relation)}.
\end{eqnarray*}
We simplify the presentation by substituting everywhere $2=3$ and $4=5$:
\begin{eqnarray*}
{}&&[2,4]=e\\
{}&&6  =  4 2  1 2^{-1} 4^{-1} \\
{}&&(1 2)^4 = (2 1)^4\\
{}&&2 = 1 2 1 2 1^{-1} 2^{-1} 1^{-1}\\
{}&&(4 6)^4 = (6 4)^4\\
{}&&4 = 6 4 6 4 6^{-1} 4^{-1} 6^{-1}\\
{}&&2^2 1 2 1 2^{-1} 1^{-1} 2^{-2} = 6 4 6  4^{-1} 6^{-1}\\
{}&&6 4^2 2^2 1 = e \ \ \mbox{(the projective relation)},
\end{eqnarray*}
and we are left with
\begin{eqnarray}
{}&&[2,4]=e\\
{}&&6 = 4 2  1 2^{-1} 4^{-1} \label{relc1}\\
{}&&(1 2)^2 = (2 1)^2\\
{}&&(4 6)^2 = (6 4)^2 \label{relc2}\\
{}&&6 4^2 2^2 1 = e.\label{relc3}
\end{eqnarray}
We substitute relation (\ref{relc1}) in relations  (\ref{relc2}) and (\ref{relc3}) and get 
$(1 4)^2 = (4 1)^2$ and $(4 2 1)^2=e$ respectively. Hence, the group $\fg{\Cc_3}$ is 
\begin{equation}
\Biggr\langle a_1, a_2, a_3  \left | 
\begin{array}{ll}
[a_2,a_3] = e    \\
(a_1 a_i)^2 = (a_i a_1)^2,  \  i=1, 2  \\
(a_1 a_2 a_3)^2 = e      
\end{array} 
\right . \Biggr\rangle, 
\end{equation}
where $a_1, a_2, a_3$ are meridians of the three quadrics in Figure \ref{C3fig}.

\begin{figure}[h]  
\epsfxsize=7cm %width  
\epsfysize=5cm %heigh  
\begin{minipage}{\textwidth}  
\begin{center}  
\epsfbox {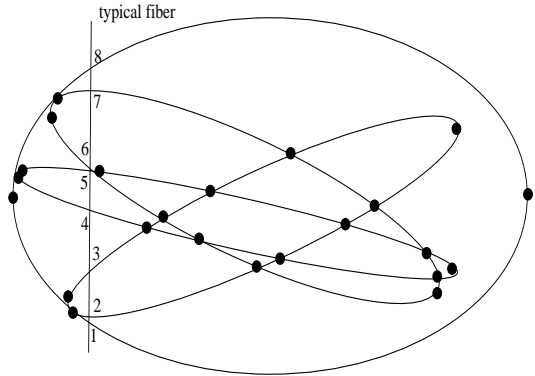}  
\end{center}  
\end{minipage}  
\caption{The arrangement $C_4$}\label{C4fig}  
\end{figure}

Let us consider now the arrangement $\Cc_4$ (Figure \ref{C4fig}).  
Following the figure, and applying the braid monodromy algorithm and the van Kampen Theorem,  
we get a presentation for  $\fg{\Cc_4}$. The initial set of generators is $1, \dots, 8$. 
The four branch points on the right side of the typical fiber contribute the relations 
$$
2=3, \ 4=5, \ 6=7 \ \ \mbox{and} \ \ 8=7 6 5 4 3 4^{-1} 6^{-1} 1 6 4 3^{-1} 4^{-1} 5^{-1} 6^{-1} 7^{-1}.
$$ 
The twelve intersections on the right side of the typical fiber contribute commutations which can be easily 
simplified to 
\begin{eqnarray*}
{}& [2,4]=[2,5]=[2,6]=[2,7]=e, \ \ [3,4]=[3,5]=[3,6]=[3,7]=e,  \\
{}& [4,6]=[4,7]=e, \ \ [5,6]=[5,7]=e.
\end{eqnarray*}
These relations enable us to simplify the relations, which we derive from the singularities on 
the left side of the typical fiber. The relations which relate to the three tangency points  are:  
$$
(1 2)^4 = (2 1)^4, \ \ (7 8)^4=(8 7)^4, \ \ (4 8)^4=(8 4)^4,
$$ 
since each tangency point is a unique common tangency between two quadrics. 
The four branch points   contribute  four relations, which can be simplified to the following forms: 
$$
(1 2)^2 = (2 1)^2, \ \ (7 8)^2=(8 7)^2, \ \ (4 8)^2=(8 4)^2, \ \ 2 1 2^{-1} = 6^{-1} 4^{-1} 8 4 6.
$$ 
The projective relation $8 7 6 5 4 3 2 1=e$ is translated to $8 {6^2} {4^2} {2^2} 1=e$.

We collect all the resulting relations, including the projective relation:
\begin{eqnarray}
{}&&8=6 4 2 1 2^{-1} 4^{-1} 6^{-1}\label{pa1}\\
{}&&[2,4]=[2,6]=[4,6]=e\\
{}&&(1 2)^2 = (2 1)^2\\
{}&&(6 8)^2=(8 6)^2\label{pa2}\\
{}&&(4 8)^2=(8 4)^2\label{pa3}\\
{}&&8 6^2 4^2 2^2 1=e. \label{pa4}
\end{eqnarray}
Now we substitute relation (\ref{pa1}) in relations (\ref{pa2}), (\ref{pa3}) and (\ref{pa4}),  
and  get $(1 6)^2 = (6 1)^2$,  $(1 4)^2 = (4 1)^2$ and $(6 4 2 1)^2=e$, respectively.  
Finally we are able to present the group $\fg{\Cc_4}$ as follows: 
\begin{equation}
\Biggr\langle a_1, a_2, a_3, a_4  \;  \left | 
\begin{array}{ll}
[a_2,a_3] = [a_2,a_4]= [a_3,a_4]= e   \\
(a_1 a_i)^2 = (a_i a_1)^2,  \ i=2, 3, 4  \\
(a_1 a_2 a_3 a_4)^2 = e                                             
\end{array} 
\right . \Biggr\rangle.
\end{equation}

In order to get the group $\fg{\Cc_n}$, we repeat the same procedure (using the Moishezon-Teicher algorithm, 
the van Kampen Theorem, and simple group simplification). We omit the lenghty computations 
and state the final result,  that $\fg{\Cc_n}$ admits the presentation (\ref{equCn}).
\hfill $\Box$

\bigskip

\noindent
\underline{\em Proof of Corollary \ref{bigCn}:}\\
The same proof as in Corollaries \ref{bigAn} and \ref{bigBn}.
\hfill $\Box$

\section{Acknowledgements}
The first author thanks the Einstein Institute, Jerusalem, and
especially Professor Ruth Lawrence-Neumark for her hospitality.
This work began while the first author was at the Mathematics
Institute, Erlangen-N\"urnberg University, Germany, and the
assistance of Professor Wolf Barth is gratefully acknowledged. 

The authors are grateful to Muhammed Uluda{\u g} for his ideas and suggestions 
concerning the types of the arrangements studied in this paper.

\end{document}